\documentclass[10pt,reqno]{amsart}
\usepackage{latexsym}
\usepackage{graphicx}
\usepackage{fancyhdr,amssymb}

\usepackage{hyperref} 

\newcommand{\N}{{\mathbb N}}

\title{Algebraic Results on SP Numbers along with a generalization}
\author{Raghavendra N Bhat, Sundarraman Madhusudanan}
\date{September 2022}

\begin{document}

\maketitle

\begin{abstract}

%-------------------------------------------------------

We defined numbers of the form $p\cdot a^2$ as \textbf{SP numbers (Square-Prime numbers)} ($a\neq1$, $p$ prime) in \cite{Bhat} along with proofs on their distribution. 
Some examples of SP Numbers : 75 = 3 $\cdot$ 25; 108 = 3 $\cdot$ 36; 45 = 5 $\cdot$ 9. 
These numbers are listed in the OEIS as A228056 \cite{OEIS}. 
In this paper, we will prove a few algebraic theorems and generalize the definition of SP Numbers to allow factors of arbitrary natural number powers.
\end{abstract}
\section{Introduction and Summary of Previous results}
In our previous paper titled `Distribution of square-prime numbers' \cite{Bhat}, we introduced the notion of a product of a prime and a square, with the square not being 1. We call these SP numbers.\\\par Here are the first 25 SP numbers:
8, 12, 18, 20, 27, 28, 32, 44, 45, 48, 50, 52, 63, 68, 72, 75, 76, 80, 92, 98, 99, 108, 112, 116, 117.\\\par In the paper we mainly proved three results. Firstly, SP numbers have an asymptotic distribution similar to prime numbers, meaning, for a large natural number $n$, the number of SP numbers smaller than $n$ is asymptotic to $\frac{n}{\log n}$, as with the case of primes. We denote the number of SP numbers smaller than $n$ by $SP(n)$. The result derived in the paper is $SP(n) = (\zeta(2) - 1)\cdot \frac{n}{\log n} + O\left(\frac{n}{\log n^2}\right).$ We see that although asymptotically similar to the distribution of primes, the actual value of $SP(n)$ is lower than $\pi(n)$ for large $n$, since $\zeta(2) - 1 < 1$.
\par Secondly, there are infinitely many pairs of SP numbers that are consecutive natural numbers. We call these SP twins. An example is (27,28). More generally, we showed that if $SP_1$ and $SP_2$ are SP Numbers, the equation $SP_1 - SP_2 = g$, where $g \in \mathbb{N}$ has infinite solutions if it has one solution. Here $g$ is the gap between the two SP Numbers.
\par Thirdly, we analyzed techniques to give asymptotic estimates for number of SP numbers having a certain last digit. For example, the number of SP numbers ending in 1 approaches $\frac{1}{400}\frac{n}{\log n}(\zeta(2,1/10) + \zeta(2,9/10) + \zeta(2,3/10) + \zeta(2,7/10) - 4).$ \par In this paper, our goal is to study algebraic properties of SP numbers, their gaps and also their structure with regard to square-free numbers. Analogies to famous prime number conjectures like the $x^2 + 1$ conjecture and primes between squares are studied. We then generalize the idea of a square multiplied by a prime to any power of a natural number multiplied by a prime.
\section{Algebraic SP Theorems}
\noindent \textbf{Theorem 2.1:} For any positive integer $x$, there exist at least one pair of SP numbers $(a,b)$ such that $a-b = x$.
\begin{proof}
We will divide all positive integers into five groups. We will then prove the existence of one pair of $SP$ numbers with difference $x$ in each group.
\begin{itemize}
    \item $\textbf{x=1:}$ 28 - 27 = 1. Thus, there exists an $SP$ pair with difference 1.\\
    \item $\textbf{x is prime: }$ Consider the following Pell's equation\\
\[m^2 - pxn^2 = 1 \tag{3}\]
where $p \in \mathbf{P} \neq x$. We know from Chakravala that $(3)$ has at least 1 positive integer solution. Let that solution be $(M,N)$. Thus,
\[M^2 - pxN^2 = 1 \tag{4}\]
\[\Rightarrow xM^2 - px^2N^2 = x \Rightarrow xM^2 - p(xN)^2 = x \tag{5}\]
Thus, we have generated two $SP$ numbers, $xM^2$ and $p(xN)^2$ that differ by $x$. Since $x$ was an arbitrary prime, this case holds.\\
\item $\textbf{x \text{is an odd composite square free number} :}$
Thus we have $x = p_1 \cdot p_2 \cdot ... \cdot p_n$, where $p_1, p_2..., p_n$ are all odd prime numbers. Let $p_1$ be the smallest prime factor of $x$. Thus $p_1 \geq 3$ and $p_2, p_3, ... p_n > 3$. Thus we can write 
\[x = p_1 \cdot y\]
where $y = p_2 \times p_3 \times ... \times p_n $. And since $p_2, p_3, ... p_n > 3$, $y > 3$. Furthermore, since $y$ is a product of odd prime numbers, $y$ has to be odd. Thus, we can also write the above equation as
\[x = p_1 \cdot (2k+1)\]
And again, since $y>3$, we have $k>1$.
\[x = p_1 \cdot [(k+1)^2 - k^2] \Rightarrow x = p_1 \cdot (k+1)^2 - p_1 \cdot k^2\]
And since $k>1$ and $p_1$ is prime, by definition $p_1 \cdot (k+1)^2$ and $p_1 \cdot k^2$ are both $SP$ numbers with difference $x$. This case holds.\\
\item $\textbf{x \text{is an even composite square free number} :}$ This means $x = p_1 \cdot p_2 \cdot ... \cdot p_n$, where $p_1 = 2$ and $p_2,p_3...p_n$ are all odd prime numbers. Thus $p_2, p_3, ... p_n \geq 3$. Thus, as previously written,
\[x = 2 \times y\]
where $y = p_2 \times p_3 \times ... \times p_n $. And since $p_2, p_3, ... p_n \geq 3$, $y \geq 3$. Furthermore, since $y$ is a product of odd prime numbers, $y$ is odd. Thus, we can also write the above equation as
\[x = 2 \cdot (2k+1)\]
And again, since $y\geq3$, we have $k\geq1$. The above equation simplifies to
\[x = 2 \cdot [(k+1)^2 - k^2] \Rightarrow x = 2 \cdot (k+1)^2 - 2 \cdot k^2\]
Note that in the previous case we had $k>1$ and here, we have $k \geq 1$. This is an issue because, when $k=1$, $2 \cdot k^2$ will be 2 which is not an $SP$ number. Note that
\[k=1 \Rightarrow y=3 \Rightarrow p_2 \cdot p_3 \cdot ... \cdot p_n = 3 \Rightarrow p_2 = 3 \]
since 3 is prime and has only 1 factor. Thus, except in the case where $p_1=2,p_2=3$, this proof works exactly like the previous case.\\
Now, let us look at $p_1=2,p_2=3$, or in other words, $x= p_1 \cdot p_2 = 6$. $18 - 12 = 6$ and $18 =2 \cdot 3^2, 12 = 3 \cdot 2^2$.\\
Since 6 was the only integer in this group for which the general method did not work, we can combine the previous 2 paragraphs to jointly say that for every even composite square free integer $x$, there exists a pair of SP numbers $(a,b)$ such that $a-b = x$.\\
\item $\textbf{x \text{is a non square free number} :}$ By the property of non square free numbers, we can write $x = t^2 \cdot s$ where $s$ is a square free number. Now, using the result from the four cases above, we can say that for all $s$, there exists a pair of $SP$ numbers $(j,k)$ such that
\[s = j - k \Rightarrow t^2 \cdot s = t^2 \cdot j - t^2 \cdot k \Rightarrow x = t^2 \cdot j - t^2 \cdot k\]
Since $(j,k)$ are $SP$ numbers, $(t^2 \cdot j$ and $t^2 \cdot k)$ are both $SP$ numbers having difference $x$. Thus, for every non square free integer $x$, there exists a pair of $SP$ numbers $(a,b)$ such that $a-b = x$.
\end{itemize}
Combining all the groups above, we can say that for any positive integer $x$, there exist at least one pair of SP numbers $(a,b)$ such that $a-b = x$.\\
\end{proof}
\noindent\textbf{Theorem 2.2:} There are infinitely many SP Numbers of the form $x^2 + 1$.
\begin{proof} Our goal is to show that the equation $x^2 +1 = pa^2$ has infinitely many solutions where $a, x \in \N$ and $p$ is prime. We start the proof with the knowledge that at least one solution already exists (we have 50, 325 and 1025 as the first 3 SP Numbers of the form $x^2 + 1$).
\[m^2 - 2n^2 = -1 \\\]
has a solution, namely $(7,5)$. Let $(A,B) = (7,5)$. Now, consider the equation
$$m^2 - 2n^2 = 1$$
We know the above Pell's equation has infinitely many solutions. Let $(X,Y)$ be an arbitrary solution to the above equation. Thus, we have the following 2 equations:
\[A^2 - 2B^2 = -1 \tag{1}\]
\[X^2 - 2Y^2 = 1 \tag{2}\]
Multiplying $(1)$ and $(2)$, we get
\[(A^2 - 2B^2)(X^2 - 2Y^2) = -1 \Rightarrow (AX + 2BY)^2 - 2(AY + BX)^2 = -1\]
Thus, for every one of the infinite solutions of the Pell's equation $(X,Y)$, we have a new solution for the equation:
\[m^2 - 2n^2 = -1 \Rightarrow 2n^2 = m^2 +1\]
Thus, there are infinitely many SP numbers of the form $x^2 +1$.\\\end{proof}
\noindent \textbf{Theorem 2.3 :} For any positive integer $x$, there exists at least one $SP$ number between $x^2$ and $(x+2)^2$.
\begin{proof}
Consider the sub-sequence of $SP$ numbers $8,18,32,50,72,98,128,162,...$. All of them are of the form $2n^2$. Notice that $2n^2$ lies between $x^2$ and $(x+2)^2$ if and only if $n \sqrt{2}$ lies between $x$ and $x+2$.\\
\indent Thus, we now need to prove that for any positive integer $x$, there exists at least one positive integer $n$ such that $n \sqrt{2}$ lies between $x$ and $(x+2)$. We will prove this using contradiction.\\
\indent Assume that there exists a positive integer $x$ such that the smallest multiple of $\sqrt{2}$ greater than $x$ is greater than $x+2$. Let this multiple be $a\sqrt{2}$ This means there is no multiple of $\sqrt{2}$ in $[x,x+2]$. But then,
\[(a-1)\sqrt{2} < x < x+2 < a\sqrt{2} \Rightarrow a\sqrt{2} - (a-1)\sqrt{2} > x+2 - x \Rightarrow \sqrt{2} > 2\]
which is a contradiction.\\
\end{proof}
\noindent \textbf{Theorem 2.4 :} All $SP$ numbers of the form $pa^2$, where $p$ is a prime and $a$ has at least one prime factor of the form $4k + 1$ can be written as a sum of 2 $SP$ numbers.
\begin{proof}
By the 2 squares theorem, if $4k+1$ is a prime, $(4k+1)^2$ can be written as $x^2 + y^2$, where $x,y > 1$. Now, since $4k+1$ is a factor of $a$, we can write
\[ a = n \cdot (4k+1) \Rightarrow pa^2 = pn^2 \cdot (4k+1)^2 \]
\[ \Rightarrow pa^2 = pn^2 \cdot (x^2+y^2) \Rightarrow pa^2 = p \cdot (nx)^2 + p  \cdot (ny)^2 \]
We now have our original $SP$ number $pa^2$ as a sum of two $SP$ numbers $p \cdot (nx)^2$ and  $p  \cdot (ny)^2$\\
\end{proof}
\section{SP Numbers of the form $x^3 + 1$}
\noindent\textbf{Conjecture 3.1}: There are infinitely many SP Numbers of the form $x^3 + 1$, if the Bunyakovsky \cite{Mathworld} conjecture is true.
\begin{proof}
Let $x = t^2-1$ for $t \in \mathbb{N}$.\\
We have $x^3 + 1 = (x+1)(x^2-x+1)\Rightarrow x^3+1 = t^2 \cdot (t^4-3t^2+3)$\\
Consider $f(t) = t^4 - 3t^2 +1$. It satisfies the following 3 properties of the Bunyakovsky conjecture:
\\1) The leading coefficient is positive.
\\2) The polynomial is irreducible over integers.
\\3) $GCD(f(2),f(3))$ is 1 and the coefficients of $f(x)$ are relatively prime.\\
\\Using the above three conditions, it follows from the Bunyakovsky conjecture that $f(t)$ is prime for infinitely many positive integers $t$. This implies $t^2 \cdot (t^4-3t^2+3)$ is of the form $prime \cdot t^2$ for infinitely many $t$ which further implies there exists infinitely many $x = t^2-1$ such that $x^3 + 1$ is an SP number.\\
\end{proof}
\section{A generalization of SP Numbers}
SP Numbers are those numbers that are of the form $pa^2$. We could extend the definition to numbers of the form $pa^k$ and get asymptotic estimates on their distribution. Let us call these numbers KP$_k$ Numbers. SP Numbers are a special case of KP numbers and can be defined to be KP$_2$. Trivially, the prime numbers and semi-prime numbers are KP$_0$ and KP$_1$ respectively. Let $KP_k(n)$ denotes the number of KP$_k$ Numbers smaller than $n$. As proved in 
\cite{Bhat}, $$KP_2(n)=(\zeta(2) - 1)\frac{n}{\log n} + O\left(\frac{n}{\log^2 n}\right)$$\\
\textbf{Theorem 4.1 :} For any $k$ greater than 1, $$KP_k(n) = (\zeta(k) - 1)\frac{n}{\log n} + O\left(\frac{n}{\log^2 n}\right)$$
\begin{proof}$$KP_k(n)=\sum_{a=2}^{\sqrt{n/2}}\pi\left(\frac{n}{a^k}\right).$$ The error term remains to be $O\left(\frac{n}{\log^2 n}\right)$, as in [1]. The main term simplifies to
$$n\sum_{a=2}^{\log^2 n}\frac{1}{a^k(\log n-\log (a^k))}=\frac{n}{\log n}\sum_{a=2}^{\log^2 n}\frac{1}{a^k(1-\frac{k\log a}{\log n})}$$
$$\Rightarrow \frac{n}{\log n}\sum_{a=2}^{\log^2 n}\frac{1}{a^k}\left(1+O\left(\frac{\log a}{\log n}\right)\right)$$
$$=\frac{n}{\log n} \sum_{a=2}^{\log^2 n}\frac{1}{a^k} + O\left(\frac{n}{\log^2 n}\right) \leq\frac{n}{\log n}\sum_{a=2}^{\infty}\frac{1}{a^k} + O\left(\frac{n}{\log^2 n}\right)$$
$$=(\zeta(k) - 1)\frac{n}{\log n} + O\left(\frac{n}{\log^2 n}\right)$$
\end{proof}
\section{PSP Numbers}
A special type of SP Numbers can be defined to be the product of a prime and a square such that the square is a prime square. An approach similar to the proof of Theorem 4.1 can be employed to get the asymptotic estimate of the count of PSP Numbers.\\\\
\textbf{Theorem 5.1 : } The number of PSP numbers smaller than $n$ is $$P(2)\frac{n}{\log n} + O\left(\frac{n}{\log^2 n}\right),$$ where $P(2)$ is the prime zeta function of 2.
\begin{proof}$$PSP(n)=\sum_{p=2}^{\sqrt{n/2}}\pi\left(\frac{n}{p^2}\right), \text{where $p$ is prime}.$$
$$\leq\frac{n}{\log n}\sum_{p=2}^{\infty}\frac{1}{p^2} + O\left(\frac{n}{\log^2 n}\right) =P(2)\frac{n}{\log n} + O\left(\frac{n}{\log^2 n}\right).$$
\end{proof}
\noindent As a corollary (owing to Theorem 4.1), we have
$$PKP_k(n)=P(k)\frac{n}{\log n} + O\left(\frac{n}{\log^2 n}\right)$$
where PKP Numbers are numbers of the form $p_1p_2^k$, where $p_1$ and $p_2$ are arbitrary prime numbers and $P(k)$ is the prime zeta function of $k$.
\section{Future work and Concluding Remarks}
\noindent Conjecture 3.1 can be explored further via the use of elliptic curves. Trying to find an SP Number of the form $x^3 + 1$ is analogous to finding integer solutions to $$x^3 + 1 = pa^2 \Rightarrow px^3 + p = (pa)^2$$
Letting $y=pa$, we have the elliptic curve equation $$y^2 = px^3 + p.$$ Computational experiments have given us over 25 integer solutions to this equation. Potential future work can involve proving the existence of infinite integer solutions to the same.
\\


\begin{thebibliography}{99}

%-------------------------------------------------------

\bibitem{Bhat}
Bhat, R. ``Distribution of Square-Prime Numbers." Missouri J. Math. Sci. 34 (1) 121 - 126, May 2022. https://doi.org/10.35834/2022/3401121 \url{https://arxiv.org/pdf/2109.10238.pdf}

\bibitem{OEIS}
The On-Line Encyclopedia of Integer Sequences, published electronically at
\url{https://oeis.org.} Sequence A228056. \url{https://oeis.org/A228056}

\bibitem{Mathworld}
Pegg, Ed Jr. ``Bouniakowsky Conjecture." From MathWorld--A Wolfram Web Resource, created by Eric W. Weisstein. \url{https://mathworld.wolfram.com/BouniakowskyConjecture.html}

\end{thebibliography}
\end{document}